\documentclass [12pt]{article}
\usepackage{amssymb,amsfonts,latexsym,amsmath}
\usepackage{comment}

\newtheorem{theorem}{Theorem}[section]%[chapter]

\newtheorem{lemma}[theorem]{Lemma}

\def\N{\mathbb{N}}

\def\Z{\mathbb{Z}}

\def\F{\mathbb{F}}

\def\cB{{\cal B}}
\def\cI{{\cal I}}
\def\cP{{\cal P}}
\def\cQ{{\cal Q}}
\def\cZ{{\cal Z}}
\def\dis{\displaystyle}
\def\sgn{{\rm sgn}}

\def\dis{\displaystyle}

\def\strictsubset{\lower .1cm \hbox{${\buildrel{\subset}\over{_{\not=
}}}$}}

\def\resp.{{\it resp.}}

\begin{document}
\title{Even perfect polynomials over $\F_2$ with four prime factors}
\author{Luis H. Gallardo - Olivier Rahavandrainy\\
Mathematics, University of Brest\\
6, Avenue Le Gorgeu, C.S. 93837,\\
29238 Brest Cedex 3, France.\\
e-mail : luisgall@univ-brest.fr - rahavand@univ-brest.fr
                                                                          \\
running head: binary perfect polynomials.\\
                                 \\
write correspondence to: Luis H. Gallardo.}
\maketitle
\newpage
%\tableofcontents
\def\cB{{\cal B}}
\def\cI{{\cal I}}
\def\cJ{{\cal J}}
\def\cP{{\cal P}}
\def\cQ{{\cal Q}}
\def\cZ{{\cal Z}}
\def\cU{{\cal U}}
\def\cX{{\cal X}}
\def\cT{{\cal T}}
\def\dis{\displaystyle}
\def\sgn{{\rm sgn}}
%\tableofcontents
\def\cB{{\cal B}}
\def\cI{{\cal I}}
\def\cJ{{\cal J}}
\def\cP{{\cal P}}
\def\cQ{{\cal Q}}
\def\cZ{{\cal Z}}
\def\cU{{\cal U}}
\def\cX{{\cal X}}
\def\cT{{\cal T}}
\def\dis{\displaystyle}
\def\sgn{{\rm sgn}}~\\
{\bf{Abstract}}\\
A perfect polynomial over the binary field  $\F_2$ is a polynomial $A \in  \F_2[x]$
that equals the sum of all its divisors. If $\gcd(A,x^2-x) \neq 1$ then we call $A$ even.
The list of  all even perfect polynomials over $\F_2$  with at most $3$ prime factors in known. The object of this paper is
to give the list of all even perfect polynomials over $\F_2$ with four prime factors.
These are all the known perfect polynomials with four prime factors over $\F_2$.

\begin{comment}
A perfect polynomial over $\F_2$ is a polynomial $A \in \F_2[x]$
that equals the sum of all its divisors. If $\gcd(A,x^2+x) \not= 1$
then we say that $A$ is even. In this paper we give the list of even
perfect polynomials with four prime divisors.
\end{comment}

{\section{Introduction}}

As usual, we denote by $\F_2$ the finite field with two elements $\{0, 1\}$.\\
For a polynomial $A \in \F_2[x]$, let $\displaystyle{\sigma(A) =
\sum_{D|A} D}$ be the sum of divisors of $A.$ We denote also, as
usual, by $\omega(A)$ the number of distinct prime (irreducible)
polynomials that divide $A.$ These two functions are multiplicative,
a fact that we shall use without more reference in the rest of the
paper. If $\sigma(A) = A,$ then we call $A$ a perfect polynomial.

\begin{comment}
The notion of perfect polynomial was introduced by Canaday
\cite{Canaday}. 
\end{comment}
The notion of perfect polynomial (over $\F_2$) was introduced by Canaday
\cite{Canaday}, the first doctoral student of Leonard Carlitz. \\ 
He studied mainly the case in which $\gcd(A,x^2+x) \neq 1.$
\begin{comment}
He studied mainly the case $\F = \F_2,$ and
$\gcd(A,x^2+x) \neq 1.$ 
\end{comment}
We may think $x^2+x \in \F_2[x]$ as being
the analogue of $2 \in \Z$ so that the ``even'' polynomials are the
polynomials with linear factors and the ``odd" ones are such that
$\gcd(A,x^2+x) =1.$ Canaday (among other results in \cite{Canaday})
classifies the even perfect polynomials with three irreducible
factors and gives without proof \cite[Theorem 11]{Canaday} the list
of all even
perfect polynomials $A$ with $\omega(A) = 4$.\\

The object of this paper (see Theorem \ref{evenwA4F2} ) is to prove Canaday's results in
\cite[Theorem 11]{Canaday}: The following polynomials are the only
even perfect polynomials $A \in \F_2[x]$ with $\omega(A)=4$ prime
factors : $$\begin{array}{l}
C_1(x) = x^2(x+1)(x^2+x+1)^2(x^4+x+1), \ C_2(x) = C_1(x+1),\\
C_3(x) =  C_3(x+1) = x^4(x+1)^4(x^4+x^3+x^2+x+1)(x^4+x^3+1),\\
C_4(x) = x^6(x+1)^3(x^3+x^2+1)(x^3+x+1), \ C_5(x) = C_4(x+1).
\end{array}$$
Observe that the two latter polynomials
\begin{comment}
Note that the two later polynomials
\end{comment}
are also perfect over $\F_4$
(see \cite{Gall-Rahav2}). 

The complete list of all even perfect polynomials over $\F_2$ with $\omega(A) \leq 4$
is then:
$$
0,1,\,\,\,  (x^2+x)^{2^{n}-1},\,\,\, T_1(x) = x^2(x+1)(x^2+x+1),\,\,\, T_1(x+1),
$$
$$
T_2(x) = x^3(x+1)^4(x^4+x^3+1),\,\,\, T_2(x+1),\,\,\,C_1(x), \ldots, C_5(x),
$$
in which $n >0$ is a positive integer.\\

In fact this list is the list of all perfect polynomials over $\F_2$ with $\omega(A) \leq 4.$
(see \cite{Gall-Rahav4}).\\

There are only two more known perfect polynomials over $\F_2$, both even,  with $\omega(A)=5$
and with degree $20$, namely:
$$
S_1(x) = x^6(x+1)^4(x^3+x+1)(x^3+x^2+1)(x^4+x^3+1),\,\,\,S_1(x+1).
$$

It may have some interest to know whether or not there are perfect polynomials over $\F_2$
with degree moderately bigger that $20$ (so that we may compute them with a computer).
These have been investigated \cite[Theorem 5.5]{Gall-Rahav3} (no solutions up to degree $28$)
 in the special case in which all exponents are equal to $2$ and the polynomial is odd.\\

{\section{Some useful facts}} 
We denote,
as usual by $\N$ the set of nonnegative integers. In this section
we recall, and we present, some necessary results for the next sections.\\

First of all, we recall some definitions and lemmata.\\
\\
{\bf{Definitions}}\\
- We define (following Canaday's terminology) as the inverse of a polynomial $P(x)$ of degree $m$, the
polynomial $\displaystyle{P^*(x)= x^m P(\frac{1}{x})}$.\\
- We say that $P$ inverts into itself if $P= P^*$.\\
- A polynomial $P$ is complete if there exists $h \in \N$ such that:
$$P = \sigma(x^h) = 1 + x + \cdots + x^h.$$

The following lemma  essentially based on a result of Dickson
(see proof of \cite[Lemma 2]{Canaday}) is key.

\begin{lemma} \label{inversion}~\\
i):  Let $P \in \F_2[x]$ be such that $P(0)=1.$ We have: $(P^*)^* = P$.\\
ii): Any complete polynomial inverts into itself.\\
iii): If $1+ x + \cdots + x^m = PQ$, where $P, Q$ are irreducible,
then either $(P=P^*, Q = Q^*)$ or $(P=Q^*, Q = P^*)$.\\
iv): If $P=P^*$, $P$ irreducible and  if $P=x^{a}(x+1)^b +1$, then:
$$P \in \{1+x+x^2, 1+x+ \cdots + x^4\}.$$
\end{lemma}
{\bf{Proof}}:\\
i) and ii) are obvious.\\
iii) follows by ii).\\
iv) is the corollary of Lemma 7 in \cite{Canaday}, (that follows from Lemma 2 of ibid.).
\begin{lemma} \label{translation}~\\
If $A(x)$ is a perfect polynomial over $\F_2$, then $A(x+1)$ is also
perfect.
\end{lemma}
\begin{lemma} {\rm (Lemma 5 in \cite{Canaday})} \label{lemma5canaday}~\\
Let $P, Q \in \F_2[x]$ and $n, m \in \N$ such that $P$ is
irreducible and $\sigma(P^{2n}) = 1 + \cdots + P^{2n} = Q^m$. Then
$m \in \{0,1\}$.
\end{lemma}
\begin{lemma} {\rm (Lemma 6 in \cite{Canaday})} \label{lemma6canaday}~\\
Let $P, Q \in \F_2[x]$ and $n, m \in \N$ such that $P$ is
irreducible and $\sigma(P^{2n}) = 1 + \cdots + P^{2n} = Q^m A, \ m>
1$. If $m$ is odd (resp. even) then  $deg(P) > (m-1) deg(Q)$ (resp.
$deg(P) > m \ deg(Q)$).
\end{lemma}
\begin{lemma} {\rm (Lemma 4 in \cite{Canaday})}
\label{lemma4canaday}~\\
If $PQ = 1+\cdots+x^{2h}$ and $P =1+ \cdots + (x+1)^{2k}$, then $h =
4$ and $k=1$; that is: $P = 1+x+x^2, \ Q = P(x^3) = 1+x^3+x^6$.
\end{lemma}
The proof of the following lemma in \cite{Canaday} uses the properties i) to iii) in Lemma \ref{inversion}: 
\begin{lemma} {\rm (Theorem 8 in \cite{Canaday})}
\label{theorem8canaday}~\\
Let $A = 1+\cdots+x^{2h} \in \F_2[x]$ such that any irreducible
factor of $A$ is of the form $x^a (x+1)^b + 1$. Then $h \in
\{1,2,3\}$.
\end{lemma}
The following crucial lemma follows from  Lemma 2.5 in \cite{Gall-Rahav2} that
says that the number of minimal primes dividing
a perfect polynomial is even:
\begin{lemma}~\\
Every even perfect polynomial $A$ over $\F_2$ with $\omega(A)=4,$ is of the form $x^h
(x+1)^k P^{l} Q^{m}$, for some odd prime polynomials $P,Q$ and for some positive integers $h,k,l,m$.
\end{lemma}

We provide proofs of the following two lemmata claimed but not proved by Canaday:

\begin{lemma} {\rm (Lemma 10 in \cite{Canaday})}
\label{lemma10canaday}~\\
Let $P \neq Q$ be two odd polynomials in $\F_2[x].$
If $x^h (x+1)^k P^{l} Q^{2n-1}$ is a perfect polynomial over $\F_2$,
and if $l \not= 2^r - 1$, then $2n-1 = 2^s - 1$.
\end{lemma}
{\bf{Proof}}:
If $l \not= 2^r - 1$ and $2n-1 \not= 2^s - 1$, then
put:
$$2n-1 = 2^s u - 1, \mbox{ where $u \geq 3$ is odd}.$$
We can write:
$$1+\cdots+Q^{2n-1} = (Q+1)^{2^s-1} (1+\cdots + Q^{u-1})^{2^s}.$$
Since $u-1 \geq 2$ is even, we have by Lemma \ref{lemma5canaday}:
$$1+\cdots + Q^{u-1} = P.$$
So, $${\rm{deg}}(Q) < {\rm{deg}}(P).$$ If $l$ is even, then by the
same argument, ${\rm{deg}}(Q) < {\rm{deg}}(P).$ It is impossible. So
$l$ is odd. We can write:
$$l = 2^r v - 1, \ 1+\cdots+P^l = (P+1)^{2^r-1}(1+\cdots +
P^{v-1})^{2^r}, \mbox{ where $v\geq 3$ is odd}.$$ Since $v-1 \geq 2$
is even, we have by Lemma \ref{lemma5canaday}:
$$1+\cdots + P^{v-1} = Q.$$
So, $${\rm{deg}}(P) < {\rm{deg}}(Q).$$ It is impossible.
\begin{flushright}
$\Box$
\end{flushright}
\begin{lemma} {\rm (Lemma 11 in \cite{Canaday})}
\label{lemma11canaday}~\\
Let $P \neq Q$ be two odd polynomials in $\F_2[x].$
If $x^h (x+1)^k P^{2l} Q^{2n-1}$ is a perfect polynomial over
$\F_2$, then $2l = 2^m$ and $m=n$.
\end{lemma}
{\bf{Proof}}:
We can write:
$$\begin{array}{l}
1+\cdots+P^{2l} = Q,\\
1+\cdots+Q^{2^n-1} = (Q+1)^{2^n-1}.
\end{array}$$
So, $P$ divides $Q+1$ and $P^2$ does not. Thus,
$$Q+1 = x^a (x+1)^b P, \mbox{ for some $a,b \in \N$}.$$
Since $\sigma(A) = A$, we obtain:
$$(1+\cdots+x^h)(1+\cdots+(x+1)^k)(x^a (x+1)^b)^{2^n-1} P^{2^n-1}Q =
x^h(x+1)^kP^{2l}Q^{2^n-1}.$$ - If $h$ and $k$ are even, then by
lemma \ref{lemma5canaday}:
$$(1+\cdots+x^h)(1+\cdots+(x+1)^k)=P^{\alpha}Q^{\beta}, \ 0 \leq \alpha,
\beta \leq 2.$$ Therefore, we must have:
$$\alpha = 1.$$
We are done.\\
\\
- If $h$ and $k$ are both odd, then by considering exponents of $P$, we see that  it
is impossible.\\
\\
- If $h$ is even and $k$ odd, then by considering exponents of $Q$, we
must have:
$$1+\cdots+x^h = P.$$
Put:
$$k+1 = 2^r u, \mbox{ where $u$ is odd}.$$
We have:
$$1+\cdots+(x+1)^k = x^{2^r-1}(1+\cdots+(x+1)^{u-1})^{2^r} =
x^{2^r-1}(P^{\gamma}Q^{\delta})^{2^r}, \ 0 \leq \gamma, \delta \leq
1.$$ - If $\gamma = 0$, then we are done.\\
\\
- If $\gamma = 1$ and $\delta = 0$, then $u-1 \geq 2$ and $n = 1$.
Thus, by considering exponents of $P$, we get:
$$l-1 = 2^{r-1}.$$
Furthermore, we can write:
$$Q+1 = P+\cdots+P^{2l} = P(1+P)(1+\cdots+P^{l-1})^2.$$
So, $l$ must be equal to $2$, and then $r = 1$, $a = 3$, $h = 4$.\\
\\
Thus:
$$P = 1+\cdots + x^4, \ Q = 1 + \cdots + P^4 = 
(1+x+x^4)(1+x+x^2+x^4+x^6+x^7+x^8+x^9+x^{12}.$$
It is impossible since $Q$ is irreducible.\\
\\
- If $\gamma = \delta = 1$, then by Lemma \ref{lemma4canaday},
${\rm{deg}}(P) = 2$ and ${\rm{deg}}(Q) = 6$. It is impossible since
$Q = 1+\cdots + P^{2l}$.\\
\\
- If $h$ is odd and $k$ even, analogous proof.
\begin{flushright}
$\Box$
\end{flushright}

In the next section we prove our main result:

\begin{theorem}
\label{evenwA4F2}
The complete list of even perfect polynomials over $\F_2$ with $4$
prime factors consists of the five polynomials $C_1(x), \ldots, C_5(x).$
\end{theorem}

\section{Perfects of the forms: $A = x^h (x+1)^k P^m Q^n$}
We may reduce (see lemmata \ref{lemma10canaday} and
\ref{lemma11canaday}) our study to the
following cases: \\

\noindent (a) $A = x^h (x+1)^k P^{2m} Q^{2n}$\\
(b) $A = x^h (x+1)^k P^{2^n} Q^{2^n-1}$\\
(c) $A = x^{2h} (x+1)^{2k} P^{2m-1} Q^{2^n-1}$\\
(d) $A = x^{2h} (x+1)^{2k-1} P^{2m-1} Q^{2^n-1}$\\
(e) $A = x^{2h-1} (x+1)^{2k-1} P^{2m-1} Q^{2^n-1}$.\\

Compare with \cite[page 733]{Canaday}.

\subsection{Case (a)}
Since $x$ and $x+1$ do not divide $\sigma(P^{2m})$, we obtain by
Lemma \ref{lemma5canaday}:
$$\sigma(P^{2m}) = 1 + \cdots + P^{2m} = Q.$$
Analogously, $$\sigma(Q^{2n}) = P.$$ Therefore, considering degrees,
we have: $$4mn = 1,$$ which is impossible.
\subsection{Case (e)}
Since $\sigma(A) = A$, we obtain:
$$x(x+1)(P+1)(Q+1) B^2 = x^{2h-1} (x+1)^{2k-1} P^{2m-1} Q^{2^n-1},
\ \mbox{for some polynomial } B.$$ It follows that $P$ (respectively
$Q$) must divide $Q+1$ (resp. $P+1$). So $P = Q+1$, which is
impossible.
\subsection{Case (b)}
We obtain:
$$\begin{array}{l}
1 + \cdots + P^{2^n} = Q, \mbox{ by Lemma \ref{lemma5canaday} }
\mbox{ and since } x, \ x+1 \mbox{ do not divide } \sigma(P^{2^n}),\\
1 + \cdots + Q^{2^n-1} =(Q+1)^{2^n-1}. \end{array}$$ Thus, $P$
divides $Q+1$ and $P^2$ does not. So, $Q$ does not divide $P+1$. We
may write:
$$\begin{array}{l}
Q+1 = P(1+P)^{2^n-1},\\
P+1 = x^{\alpha}(x+1)^{\beta}, \ \alpha, \beta \geq 1 \end{array}$$
\subsubsection{Case $h$, $k$ even} The two monomials $x$ and $x+1$ do
not divide $\sigma(x^h), \sigma((x+1)^k)$. So:
$$\begin{array}{l}
1+ \cdots + x^{h} = P^{a_0}Q^{b_0}, \ a_0, b_0 \in \{0,1\},\\
1+ \cdots + (x+1)^{k} = P^{a_1}Q^{b_1}, \ a_1, b_1 \in \{0,1\}.
\end{array}$$
Since $\sigma(A) = A$, we obtain:
$$\begin{array}{l}
Q(Q+1)^{2^n-1}P^{l}Q^{r} = x^h (x+1)^k P^{2^n} Q^{2^n-1},\\
l = a_0+a_1, \ r = b_0+b_1, \ l, r \in \{0, 1, 2\}.
\end{array}$$
Considering the exponents of $P$ and $Q$, we have: $$2^n-1 + l =
2^n, \ r+1 = 2^n - 1.$$ So, $$l = 1, \ n \in \{1,2\}.$$
(i)- \underline{Case $n = 1$}:\\
We have: $$r = 0, \ 1+P+P^2 = Q, \ 1+ \cdots + x^{h} = P = 1+ \cdots
+ (x+1)^{k}, \ h = k.$$ Since $P = x^{\alpha}(x+1)^{\beta} +1$, by
Lemma \ref{theorem8canaday}, $P \in \{1+x+x^2, 1+
\cdots+x^4\}$.\\
- If $P= 1+x+x^2$, then $h=k=2, \ Q = 1+x+x^4$. Thus $A =
x^2(x+1)^2P^2Q$ which is not perfect.\\
- If $P= 1+\cdots +x^4$, then $Q = 1+P+P^2 =
(1+x+x^2)(1+x^2+x^4+x^5+x^6)$ is reducible. It is impossible.\\
\\
(ii)- \underline{Case $n = 2$}:\\
We have: $$\begin{array}{l} r = 2, \ 1+\cdots+P^4 = Q,\\
 1+ \cdots +
x^{h} = PQ,\\
1+ \cdots + (x+1)^{k} = Q. \end{array}$$ By Lemma
\ref{lemma4canaday}, we have:
$$h = 8, \ k = 2, \ Q = 1+x+x^2, \ P = 1 + x^3 +x^6.$$
So, $Q \not= 1+\cdots + P^4$. It is impossible.

\subsubsection{Case $h$, $k$ odd}
Since $\sigma(A) = A$, we obtain:
\begin{equation}
\label{hkodd}
x(x+1)Q(Q+1)^{2^n-1}B^2 = x^h (x+1)^k P^{2^n} Q^{2^n-1}.
\end{equation}
Since $P$ divides $Q+1$ and $P^2$ does not, by considering the exponent
of $P$, we see that the equality (\ref{hkodd}) is impossible.
\subsubsection{Case $h$ odd, $k$ even}
Put $h = 2l-1$ and $k = 2r$. \\
By Lemma \ref{lemma6canaday}, we have:
$$1+\cdots+(x+1)^k = 1 + \cdots + (x+1)^{2r} = P^aQ^b,
\mbox{ for some } a,b\in \{0,1\}.$$ Since $\sigma(A) = A$, we
obtain:
$$(x+1)(1+\cdots+x^{l-1})^2Q(Q+1)^{2^n-1}P^aQ^b =
x^{2l-1} (x+1)^{2r} P^{2^n} Q^{2^n-1}.$$ Since $P$ divides $Q+1$ and
$P^2$ does not, if $b = 1$ (resp. $a = 0$), then the exponent of $Q$
(resp. of $P$) in the right hand side is even (resp. odd). It is
impossible. So, $b = 0$ and $a = 1$. Therefore: $$P=1+\cdots
+(x+1)^{2r} = x^{\alpha}(x+1)^{\beta} +1.$$
By Lemma \ref{theorem8canaday}, $P \in \{1+x+x^2, 1+x^3+x^4\}$.\\
\\
(i)- \underline{Case $P = 1+x+x^2$}:\\
We have $k = 2r = 2$, and by considering the exponent of $x+1$ we get:
$$n=1, \ Q = 1+P+P^2 = 1+x+x^4.$$
So, $$l=1, \ h = 1.$$ We obtain the polynomial $C_1(x)$, and by
Lemma \ref{translation}, we get the polynomial $C_1(x+1)$.\\
\\
(ii)- \underline{Case $P = 1+x^3+x^4$}:\\
We have: $$2r = 4, \ Q+1 = (1+P)^{2^n-1}P =
x^{3(2^n-1)}(x+1)^{2^n-1}P.$$ By considering the exponent of $x+1$, we
have:
$$(2^n-1)^2 +1 \leq 4$$
So, $$n = 1,$$ and $$Q = 1+P+P^2 = 1+x^3+x^4+x^6+x^8
=(1+x+x^2)(1+x+x^4+x^5+x^6).$$ It is impossible.

\subsection{Case (c)}
By Lemma \ref{lemma6canaday}, we obtain:
$$\begin{array}{l}
1 + \cdots + x^{2h} = P^{a_0}Q^{b_0}, \\
1 + \cdots + (x+1)^{2k} = P^{a_1}Q^{b_1},\\
a_0, b_0, a_1, b_1 \in \{0,1\}. \end{array}$$ Since $\sigma(A) = A$,
we obtain:
$$(P+1)(Q+1)^{2^n-1}P^{a_0+a_1}Q^{b_0+b_1}(1 + \cdots + P^{m-1})^2 =
x^{2h} (x+1)^{2k} P^{2m-1} Q^{2^n-1}.$$ Thus:
$$\begin{array}{l}
1+P =x^{\alpha_1} (x+1)^{\beta_1}Q^{\gamma_1},\\
1+Q =x^{\alpha_2} (x+1)^{\beta_2}P^{\gamma_2},\\
\alpha_1, \beta_1, \gamma_1, \alpha_2, \beta_2, \gamma_2 \in \N.
\end{array}$$ We can reduce the work to three cases, since the integers $h$
and $k$ play symmetric roles (by Lemma \ref{translation}).

\subsubsection{Case $a_0 = b_0 = b_1 = 1, a_1 = 0$}
We have:
$$\begin{array}{l}
1 + \cdots + x^{2h} = PQ, \\
1 + \cdots + (x+1)^{2k} = Q. \end{array}$$ So, by Lemma
\ref{lemma4canaday}, we obtain:
$$Q = 1+x+x^2, \ P = 1+x^3+x^6, \ h = 4, \ k = 1.$$
Since $\sigma(A) = A$, by considering the exponent of $x+1$,
we obtain:
$$n = 1,$$
and thus:
$$x^4(x+1)^2PQ^3(1+\cdots+P^{m-1})^2 = x^{8} (x+1)^{2} P^{2m-1}
Q.$$ Thus, $x$ must divide $B = 1+\cdots+P^{m-1}$. So, $x+1$ must
divide $B$. It is impossible.

\subsubsection{Case $a_0 = b_0 = a_1 = 1, b_1 = 0$}
We have:
$$\begin{array}{l}
1 + \cdots + x^{2h} = PQ, \\
1 + \cdots + (x+1)^{2k} = P. \end{array}$$ So, by Lemma
\ref{lemma4canaday}, we obtain:
$$P = 1+x+x^2, \ Q = 1+x^3+x^6, \ h = 4, \ k = 1.$$
We obtain the same contradiction as in the previous case.

\subsubsection{Case $a_0 = b_1 = 1, a_1 = b_0 = 0$}
We have:
$$\begin{array}{l}
1 + \cdots + x^{2h} = P, \\
1 + \cdots + (x+1)^{2k} = Q \end{array}$$~\\

Therefore, the monomials $x, \ x+1$ divide $P+1$ and $Q+1$. But
$x^2$ (resp. $(x+1)^2$) does not divide $P+1$ (resp. $Q+1$).\\
Since $\sigma(A) = A$, we have:
\begin{equation} \label{casc1}
PQ(P+1)(Q+1)^{2^n-1}(1+\cdots+P^{m-1})^2 = x^{2h} (x+1)^{2k}
P^{2m-1} Q^{2^n-1}
\end{equation}
(i)- \underline{Case $\gamma_1 = \gamma_2 = 0$}:\\
\\
In this case, $P$ does
not divide $Q+1$ and $Q$ does not divide $P+1$.\\
We obtain: $$P+1 = x (x+1)^{\beta_1}, \ Q+1 = x^{\alpha_2} (x+1).$$
Therefore, by relation (\ref{casc1}):
$$m = 1 \mbox{ and } n = 1.$$ So, by Lemma \ref{theorem8canaday}:
$$P \in \{1+x+x^2, 1+\cdots +x^4\}, \ Q \in \{1+x+x^2,
1+x^3+x^4\}.$$
We must have:
$$P = 1 + \cdots +x^4, \ Q = 1+x^3+x^4.$$ So, $$h = k =
2.$$ We get
the polynomial $C_3(x)$, and thus the polynomial $C_3(x+1) = C_3(x)$.\\
\\
(ii)- \underline{Case $\gamma_1 = 0, \ \gamma_2 \geq 1$}:\\
\\
The polynomial $P$ divides $Q+1$, and by relation (\ref{casc1}), the
integer $\gamma_2$ must be even. So:
$$Q+1 = x^{\alpha_2} (x+1) P^{2u}, \ u \geq 1.$$ In particular,
$P^2$ divides $Q+1$.\\
\\
Furthermore, $Q$ does not divide $P+1$. So,
$$P=x(x+1)^{\beta_1} + 1.$$
So, by Lemma \ref{theorem8canaday}, $P \in
\{1+x+x^2, 1+\cdots +x^4\}$.\\
\\
- If $P= 1+x+x^2$, then $2h = 2, \ n=1$ and $\alpha_2 = 1$ (consider the
exponents of $x$ in the relation (\ref{casc1})). We can write:
$$Q+1 = x (x+1)P^{2u}.$$
By considering the exponent of $P$, we have: $$u = m-1.$$ and thus:
$$m \geq 2.$$
Moreover, the relation (\ref{casc1}) becomes:
$$(1 + \cdots +
P^{m-1})^2 = (x+1)^{2k-2}.$$ So,
$$k = 1, \ m = 1.$$ It is impossible.\\
\\
- If $P= 1+\cdots+x^4$, then $2h = 4$. \\
We can write:
$$Q+1 = x^{\alpha_2}(x+1)P^{2u},
\mbox{ where $\alpha_2$ is odd and $u \geq 1$}.$$ By considering the
exponent of $x$ in relation (\ref{casc1}), we have:
$$\mbox{ either $(n = 2, \ \alpha_2 = 1)$~or~$(n=1,~
\alpha_2~\in~\{1,3\})$.}$$
\underline{Case $n=2, \ \alpha_2 = 1$}:\\
By considering the exponent of $P$, we have:
$$m = 3u+1 \geq 4.$$
Moreover, we must have:
$$1+\cdots + P^{m-1} = (x+1)^{k-3}Q.$$
So, $$k = 3, \ 1+\cdots + P^{m-1} = Q.$$ Thus, $P^2$ does not divide
$Q+1$. It is impossible.\\
\\
\underline{Case $n=1$}:\\
By considering the exponent of $P$, we have:
$$u = m-1 \geq 1.$$
Moreover, we must have:
$$(1+\cdots + P^{m-1})^2 = x^{4-\alpha_2-1}(x+1)^{2k-4}.$$
Thus, $m-1$ is odd, $\alpha_2 = 1$. \\
By writing: $$1+\cdots + P^{m-1} = (1+P)(1+\cdots + P^{m/2-1})^2.$$ We
must have: $m=2$, $u=1$, and $k=5$. So,
$$\begin{array}{l}
Q = 1 + \cdots + (x+1)^{10} = 1+x+x^2+x^7+x^8+x^9+x^{10},\\
Q+1 = x(x+1)P^2 = x(x+1)(1+\cdots + x^4)^2 = x + \cdots + x^{10}.
\end{array}$$
It is impossible.\\
\\
(iii)- \underline{Case $\gamma_1 \geq 1, \ \gamma_2 = 0$}:\\
\\
In this case, $P$ does not divide $Q+1$, and $Q$ divides
$P+1$.\\
So,
$$Q=x^{\alpha_2}(x+1) + 1.$$
So, by Lemma \ref{theorem8canaday}:
$$Q \in
\{1+x+x^2, 1+ x^3 +x^4\}.$$ Therefore, by relation (\ref{casc1}):
$$m = 1.$$
So:
$$(P+1)(Q+1)^{2^n-1}=
x^{2h} (x+1)^{2k}Q^{2^n-2}.$$ - If $Q= 1+x+x^2$, then $k=1$,
$\gamma_1 = 2u = 2^n - 2$ is even, and $\beta_1$ is odd. We can
write:
$$P+1 = x(x+1)^{\beta_1} Q^{\gamma_1}, \ \gamma_1 = 2u \geq 2.$$
Considering the exponent of $x+1$, we have:
$$2 = 2^{n}-1 + \beta_1.$$ So:
$$n= \beta_1 = 1.$$
Thus:
$$\gamma_1 = 2^n - 2 = 0.$$ It is impossible.\\
\\
- If $Q= 1+x^3+x^4$, then $$2k=4 = 2^{n}-1 + \beta_1.$$ So:
$$\mbox{ either $(n=1, \ \beta_1 = 3)$ or $(n=2, \ \beta_1 = 1)$}.$$
The first case is impossible since $\gamma_1 = 2^n - 2 \geq 2$.\\
So, $$n=2, \ \beta_1 = 1, \gamma_1 = 2, \ 2h = 3.(2^n - 1) + 1 =
10.$$ Thus:
$$\begin{array}{l}
P = 1 + \cdots + x^{10},\\
P+1 = x(x+1)Q^2 = x(x+1)(1+ x^3 + x^4)^2 = x+x^2+x^7+x^8+x^9+x^{10}.
\end{array}$$
It is impossible.
\subsection{Case (d)}
We obtain:
$$\begin{array}{l}
1 + \cdots + x^{2h} = P^{a_0}Q^{b_0}, \ a_0, b_0 \in \{0,1\}
\mbox{
by
Lemma \ref{lemma6canaday}},\\
1 + \cdots + (x+1)^{2k-1} = x(1+\cdots+(x+1)^{k-1})^2.
\end{array}$$ Since $\sigma(A) = A$, we obtain:
\begin{equation} \label{casd1}
x(P+1)(Q+1)^{2^n-1}B^2P^{a_0}Q^{b_0} = x^{2h} (x+1)^{2k-1} P^{2m-1}
Q^{2^n-1} \end{equation} Thus:
$$\begin{array}{l}
1+P =x^{\alpha_1} (x+1)^{\beta_1}Q^{\gamma_1},\\
1+Q =x^{\alpha_2} (x+1)^{\beta_2}P^{\gamma_2}. \end{array}$$
By considering degrees, we obtain: $$\gamma_1 \gamma_2 \leq 1.$$
If $\gamma_1 = \gamma_2 = 1$, then $Q =  P+1$. It is impossible.\\
So, $\gamma_1 \gamma_2 = 0$. We have three cases:

\subsubsection{Case: $\gamma_1 = \gamma_2 = 0$}
In this case, $Q$ (resp. $P$) does not divide $P+1$ (resp. $Q+1$).
We may write: $$P = x^{\alpha_1} (x+1)^{\beta_1} + 1, \ Q =
x^{\alpha_2} (x+1)^{\beta_2} + 1.$$ - If $1 + \cdots + x^{2h} = P$,
then the relation (\ref{casd1}) becomes:
$$x(P+1)(Q+1)^{2^n-1}B^2P =
x^{2h} (x+1)^{2k-1} P^{2m-1} Q^{2^n-1}.$$ It is impossible (consider
the exponent of $Q$).\\
- If $1 + \cdots + x^{2h} = Q$, then the relation (\ref{casd1})
becomes:
$$x(P+1)(Q+1)^{2^n-1}B^2Q =
x^{2h} (x+1)^{2k-1} P^{2m-1} Q^{2^n-1}.$$ It is impossible (consider
the exponent of $P$).\\
- If $1 + \cdots + x^{2h} = PQ$, then by Lemma
\ref{theorem8canaday}:
$$P, Q \in \{x^3 +x^2 +1, x^3+x+1\}, \ h = 3.$$
We get the polynomial $C_4(x)$ and thus also the polynomial $C_5(x) = C_4(x+1)$.

\subsubsection{Case: $\gamma_1 = 0, \ \gamma_2 \geq 1$}
In this case, we may write: $$P = x^{\alpha_1} (x+1)^{\beta_1} + 1,
\ Q = x^{\alpha_2} (x+1)^{\beta_2}P^{\gamma_2} + 1.$$ So
${\rm{deg}}(P) < {\rm{deg}}(Q)$.\\
\\
- If $1 + \cdots + x^{2h} = P$, then it is impossible as in the above
case (consider the exponent of $Q$).\\
- If $1 + \cdots + x^{2h} = Q$, then:
$$\begin{array}{l}
a_0 = 1, \ b_0 = 0,\\
\mbox{ $x$ divides $Q+1$, $x^2$ does not},\\
Q+1 = x(x+1)(1+ \cdots + x^{h-1})^2.
\end{array}$$
So, $\alpha_2 = 1$ and $\gamma_2$ is even.\\
By considering the exponent of $P$, we see that the relation
(\ref{casd1}) does not hold. It is impossible.\\
- If $1 + \cdots + x^{2h} = PQ$, then by Lemma \ref{inversion},
since ${\rm{deg}}(P) < {\rm{deg}}(Q)$, the polynomial $P$ (resp.
$Q$) inverts into itself, and $P \in \{1+x+x^2, 1+
\cdots + x^4\}$.\\
Therefore, $\alpha_1 = 1$ and $\beta_1 \in\{1,3\}$ is odd. Thus,
by considering the equality:
$$x(P+1)(Q+1)^{2^n-1}B^2PQ =
x^{2h} (x+1)^{2k-1} P^{2m-1} Q^{2^n-1},$$ we obtain that the integers $\alpha_2,
\beta_2$ and $\gamma_2$ must be even. So, $Q+1$ is a square. It is
impossible by the irreducibility of $Q$.
\subsubsection{Case: $\gamma_1 \geq 1, \ \gamma_2 = 0$}
In this case, we may write: $$P = x^{\alpha_1}
(x+1)^{\beta_1}P^{\gamma_1} + 1, \ Q = x^{\alpha_2} (x+1)^{\beta_2}
+ 1.$$ The proof is analogous to that of the previous case, by switching $P$ and $Q$.\\
\\
- If $1 + \cdots + x^{2h} = Q$, then it is impossible (consider the exponent of $P$).\\
- If $1 + \cdots + x^{2h} = P$, then:
$$\begin{array}{l}
\mbox{ $x$ divides $P+1$, $x^2$ does not},\\
P+1 = x(x+1)(1+ \cdots + x^{h-1})^2.
\end{array}$$
So, $\alpha_1 = 1$ and $\gamma_1$ is even.\\
By considering the exponent of $Q$, we see that the following equality
does not hold:
$$x(P+1)(Q+1)^{2^n-1}B^2P =
x^{2h} (x+1)^{2k-1} P^{2m-1} Q^{2^n-1}.$$ It is impossible.\\
\\
- If $1 + \cdots + x^{2h} = PQ$, then by Lemma \ref{inversion},
since ${\rm{deg}}(Q) < {\rm{deg}}(P)$, the polynomial $P$ (resp.
$Q$) inverts into itself, and $Q \in \{1+x+x^2, 1+x+
\cdots + x^4\}$.\\
Therefore, $\alpha_2 = 1$ and $\beta_2 \in\{1,3\}$ is odd. Thus,
by considering the equality:
$$x(P+1)(Q+1)^{2^n-1}B^2PQ =
x^{2h} (x+1)^{2k-1} P^{2m-1} Q^{2^n-1},$$ the integers $\alpha_1,
\beta_1$ and $\gamma_1$ must be even. So, $P+1$ is a square. It is
impossible by the irreducibility of $P$.\\

This finishes the proof of Theorem \ref{evenwA4F2}.

%\newpage
\def\thebibliography#1{\section*{\titrebibliographie}
\addcontentsline{toc}
{section}{\titrebibliographie}\list{[\arabic{enumi}]}{\settowidth
 \labelwidth{[
#1]}\leftmargin\labelwidth \advance\leftmargin\labelsep
\usecounter{enumi}}
\def\newblock{\hskip .11em plus .33em minus -.07em} \sloppy
\sfcode`\.=1000\relax}
\let\endthebibliography=\endlist

\def\biblio{\def\titrebibliographie{References}\thebibliography}
\let\endbiblio=\endthebibliography

%\def\references{\def\titrebibliographie{R\' ef\'
%erences}\thebibliography}
%\let\endreferences=\endthebibliography

%%%% MACROS DE SEROUL POUR LES REFERENCES %%%%

%%%%%%% bibliographie selon AMS style %%%%%%%%%%%
%%%%%%% inspir de TUGboat 11 (1990), p. 609 %%%%%%%

\newbox\auteurbox
\newbox\titrebox
\newbox\titrelbox
\newbox\editeurbox
\newbox\anneebox
\newbox\anneelbox
\newbox\journalbox
\newbox\volumebox
\newbox\pagesbox
\newbox\diversbox
\newbox\collectionbox
%--------------------------------------------
\def\fabriquebox#1#2{\par\egroup
\setbox#1=\vbox\bgroup \leftskip=0pt \hsize=\maxdimen \noindent#2}
%--------------------------------------------
\def\bibref#1{\bibitem{#1}

\mbox{}\ignorespaces

\setbox0=\vbox\bgroup}
%--------------------------------------------
\def\auteur{\fabriquebox\auteurbox\styleauteur}
\def\titre{\fabriquebox\titrebox\styletitre}
\def\titrelivre{\fabriquebox\titrelbox\styletitrelivre}
\def\editeur{\fabriquebox\editeurbox\styleediteur}

\def\journal{\fabriquebox\journalbox\stylejournal}

\def\volume{\fabriquebox\volumebox\stylevolume}
\def\collection{\fabriquebox\collectionbox\stylecollection}
%--------------------------------------------
{\catcode`\- =\active\gdef\annee{\fabriquebox\anneebox\catcode`\-
=\active\def -{\hbox{\rm
\string-\string-}}\styleannee\ignorespaces}}
%--------------------------------------------
{\catcode`\-
=\active\gdef\anneelivre{\fabriquebox\anneelbox\catcode`\-=
\active\def-{\hbox{\rm \string-\string-}}\styleanneelivre}}
%--------------------------------------------
{\catcode`\-=\active\gdef\pages{\fabriquebox\pagesbox\catcode`\-
=\active\def -{\hbox{\rm\string-\string-}}\stylepages}}
%--------------------------------------------
{\catcode`\-
=\active\gdef\divers{\fabriquebox\diversbox\catcode`\-=\active
\def-{\hbox{\rm\string-\string-}}\rm}}
%--------------------------------------------
\def\ajoutref#1{\setbox0=\vbox{\unvbox#1\global\setbox1=
\lastbox}\unhbox1 \unskip\unskip\unpenalty}
%--------------------------------------------
\newif\ifpreviousitem
\global\previousitemfalse
\def\separateur{\ifpreviousitem {,\ }\fi}
%--------------------------------------------
\def\voidallboxes
{\setbox0=\box\auteurbox \setbox0=\box\titrebox
\setbox0=\box\titrelbox \setbox0=\box\editeurbox
\setbox0=\box\anneebox \setbox0=\box\anneelbox
\setbox0=\box\journalbox \setbox0=\box\volumebox
\setbox0=\box\pagesbox \setbox0=\box\diversbox
\setbox0=\box\collectionbox \setbox0=\null}
%--------------------------------------------
\def\fabriquelivre
{\ifdim\ht\auteurbox>0pt
\ajoutref\auteurbox\global\previousitemtrue\fi
\ifdim\ht\titrelbox>0pt
\separateur\ajoutref\titrelbox\global\previousitemtrue\fi
\ifdim\ht\collectionbox>0pt
\separateur\ajoutref\collectionbox\global\previousitemtrue\fi
\ifdim\ht\editeurbox>0pt
\separateur\ajoutref\editeurbox\global\previousitemtrue\fi
\ifdim\ht\anneelbox>0pt \separateur \ajoutref\anneelbox
\fi\global\previousitemfalse}
%--------------------------------------------
\def\fabriquearticle
{\ifdim\ht\auteurbox>0pt        \ajoutref\auteurbox
\global\previousitemtrue\fi \ifdim\ht\titrebox>0pt
\separateur\ajoutref\titrebox\global\previousitemtrue\fi
\ifdim\ht\titrelbox>0pt \separateur{\rm in}\
\ajoutref\titrelbox\global \previousitemtrue\fi
\ifdim\ht\journalbox>0pt \separateur
\ajoutref\journalbox\global\previousitemtrue\fi
\ifdim\ht\volumebox>0pt \ \ajoutref\volumebox\fi
\ifdim\ht\anneebox>0pt  \ {\rm(}\ajoutref\anneebox \rm)\fi
\ifdim\ht\pagesbox>0pt
\separateur\ajoutref\pagesbox\fi\global\previousitemfalse}
%--------------------------------------------
\def\fabriquedivers
{\ifdim\ht\auteurbox>0pt
\ajoutref\auteurbox\global\previousitemtrue\fi
\ifdim\ht\diversbox>0pt \separateur\ajoutref\diversbox\fi}
%--------------------------------------------
\def\endbibref
{\egroup \ifdim\ht\journalbox>0pt \fabriquearticle
\else\ifdim\ht\editeurbox>0pt \fabriquelivre
\else\ifdim\ht\diversbox>0pt \fabriquedivers \fi\fi\fi
.\voidallboxes}
%--------------------------------------------

\let\styleauteur=\sc
\let\styletitre=\it
\let\styletitrelivre=\sl
\let\stylejournal=\rm
\let\stylevolume=\bf
\let\styleannee=\rm
\let\stylepages=\rm
\let\stylecollection=\rm
\let\styleediteur=\rm
\let\styleanneelivre=\rm

\begin{biblio}{99}
%\addcontentsline{toc}{section}{Bibliographie}

\begin{bibref}{Canaday}
\auteur{E. F. Canaday} \titre{The sum of the divisors of a
polynomial} \journal{Duke Math. Journal} \volume{8} \pages 721 - 737
\annee 1941
\end{bibref}

\begin{bibref}{Beard}
\auteur{T. B. Beard Jr, James. R. Oconnell Jr, Karen I. West}
\titre{Perfect polynomials over $GF(q)$} \journal{Rend. Accad.
Lincei} \volume{62} \pages 283 - 291 \annee 1977
\end{bibref}

\begin{bibref}{Gall-Rahav}
\auteur{L. Gallardo, O. Rahavandrainy} \titre{On perfect polynomials
over $\F_4$} \journal{Portugaliae Mathematica} \volume{62 - Fasc. 1}
\pages 109- 122 \annee 2005
\end{bibref}

\begin{bibref}{Gall-Rahav2}
\auteur{L. Gallardo, O. Rahavandrainy} \titre{Perfect polynomials
over $\F_4$ with less than five prime factors} \journal{Portugaliae 
Mathematica}
\volume{64 - Fasc. 1} \pages 21-38 \annee 2007
\end{bibref}

\begin{comment}
\begin{bibref}{Gall-Rahav2}
\auteur{L. Gallardo, O. Rahavandrainy} \titre{Perfect polynomials
over $\F_4$ with less than five prime factors} \journal{to appear in
Portugaliae Mathematica} \volume{} \pages \annee 2006
\end{bibref}
\end{comment}

\begin{bibref}{Gall-Rahav3}
\auteur{L. H. Gallardo, O. Rahavandrainy} \titre{Odd perfect polynomials
over $\F_2$} \journal{J.T.N.B.}
\volume{19} \pages 167-176 \annee 2007
\end{bibref}

\begin{bibref}{Gall-Rahav4}
\auteur{L. H. Gallardo, O. Rahavandrainy} \titre{There is no odd perfect polynomial
over $\F_2$ with four prime factors} \journal{Preprint}
%\volume{19} \pages 167-176 \annee 2007
\end{bibref}

\begin{comment}
\begin{bibref}{Rudolf}
\auteur{Rudolf Lidl, Harald Niederreiter} \titrelivre{Finite Fields,
Encyclopedia of Mathematics and its applications} \editeur{Cambridge
University Press} \anneelivre 1983 (Reprinted 1987)
\end{bibref}

\begin{bibref}{Steuerw}
\auteur{Rudolf Steuerwald} \titre{Versch\"arfung einer notwendigen
Bedingung f\"ur die Existenz einer ungeraden vollkommenen Zahl}
\journal{S. B. math.-nat. Abt. Bayer. Akad. Wiss. M\"unchen}
\volume{} \pages 69 - 72 \annee 1937
\end{bibref}
\end{comment}

\end{biblio}

\end{document}